\documentclass[12pt]{article}
\usepackage{amsmath}
\usepackage{amsthm}
\usepackage{amsfonts}
\usepackage{amssymb}
\usepackage{esint}
\usepackage{eucal}
\usepackage{mathrsfs}
\usepackage{graphicx,graphics}
\numberwithin{equation}{section}
\usepackage{graphicx,graphics}
\usepackage{pdfsync}
\usepackage{multirow}
\usepackage{endnotes}
\def \dis {\displaystyle}

\def \limi {\;\mathop{\longrightarrow}_{n\to\infty}\;}

\def \confai {-\kern -.5em\rightharpoonup}
\def \cqfd {\hfill$\Box$}
\def \div{\mbox{\rm div}}

\def \ga {\gamma}
\def \Ga {\Gamma}

\def \De {\Delta}

\def \Om {\Omega}
\def \la {\lambda}

\def \ph {\varphi}

\def \si {\sigma}

\def \NN {\mathbb N}

\def \PP {\mathbb P}

\def \RR {\mathbb R}

\def \D {\mathscr{D}}

\def \beq {\begin{equation}}
\def \eeq {\end{equation}}
\def \ba {\begin{array}}
\def \ea {\end{array}}
\def \bs {\bigskip}
\def \ms {\medskip}

\def \ecart {\noalign{\medskip}}
\newtheorem{Thm}{Theorem}[section]

\newtheorem{Pro}[Thm]{Proposition}
\newtheorem{Lem}[Thm]{Lemma}

\newtheorem{Adef}[Thm]{Definition}
\newenvironment{Def}{\begin{Adef}\rm}{\end{Adef}}
\newtheorem{Arem}[Thm]{Remark}
\newenvironment{Rem}{\begin{Arem}\rm}{\end{Arem}}
\newtheorem{Aexa}[Thm]{Example}
\newenvironment{Exa}{\begin{Aexa}\rm}{\end{Aexa}}
\newtheorem{Anot}[Thm]{Notation}

\def \refe #1.{(\ref{#1})}
\def \reff #1.{figure~\ref{#1}}
\def \refs #1.{Section~\ref{#1}}
\def \refss #1.{Subsection~\ref{#1}}
\def \refD #1.{Definition~\ref{#1}}
\def \refT #1.{Theorem~\ref{#1}}
\def \refL #1.{Lemma~\ref{#1}}
\def \refC #1.{Corollary~\ref{#1}}
\def \refP #1.{Proposition~\ref{#1}}
\def \refPt #1.{Properties~\ref{#1}}
\def \refR #1.{Remark~\ref{#1}}
\def \refE #1.{Example~\ref{#1}}
\def \refN #1.{Notation~\ref{#1}}
\newcounter{marnote}

\usepackage[colorlinks=true, linkcolor=blue]{hyperref} 

\title{Reconstruction of isotropic conductivities \\ from non smooth electric fields}
\author{Marc Briane
\\
\normalsize Univ Rennes, INSA Rennes,  CNRS, IRMAR - UMR 6625, F-35000 Rennes, France
\\
\normalsize mbriane@insa-rennes.fr
}
\begin{document}
\maketitle
\begin{abstract}
In this paper we study the isotropic realizability of a given non smooth gradient field $\nabla u$ defined in $\RR^d$, namely when one can reconstruct an isotropic conductivity $\si$ such that $\si\nabla u$ is divergence free in $\RR^d$. On the one hand, in the case where $\nabla u$ is non-vanishing, uniformly continuous in $\RR^d$ and $\De u$ is a bounded function in $\RR^d$, we prove the isotropic realizability of $\nabla u$ using the associated gradient flow combined with the DiPerna, Lions approach for solving ordinary differential equations in suitable Sobolev spaces. On the other hand, in the case where the gradient $\nabla u$ is piecewise regular, we prove roughly speaking that the isotropic realizability holds if and only if the normal derivatives of $u$ on each side of the gradient discontinuity interfaces have the same sign. Some examples of conductivity reconstruction are given.
\end{abstract}
\par\bs\noindent
\section{Introduction}
In Electrophysics there are some constraints implicitly satisfied by the electric field in a prescribed conductive material.
For example, Alessandrini and Nesi \cite{AlNe} have shown that a smooth periodic electric field cannot vanish in dimension two, while it may vanish in dimension three as proved in \cite{Anc,BMN}. This three-dimensional specificity of the electric field allows us to derive a surprising property of the Hall effect: the sign of the Hall voltage is indeed inverted in a three-dimensional {\em metamaterial} inspired by a chain mail armor. The anomalous Hall effect has been first proved theoretically in~\cite{BrMi}, then it has been simplified and validated experimentally in~\cite{KKW}. Very recently it has been emphasized simultaneously in {\em Physics Today}~\cite{Mil} and {\em Nature}~\cite{Not}.
\par
Conversely, {starting from a regular gradient field $\nabla u\neq 0$ in $\RR^d$
(\footnote{When $d=2$, $\nabla u\neq 0$ in the periodic case (see \cite{AlNe}), otherwise it is obvious that there exist solutions with $\nabla u$ vanishing somewhere. A treatment of such cases can be found in \cite{Ale}. The case $d=3$ is quite different, since $\nabla u$ may vanish somewhere in the periodic case (see \cite{BMN}).})} the natural inverse problem is to reconstruct from $\nabla u$ a possibly isotropic conductivity $\sigma$ which satisfies the conductivity equation
\beq\label{conequ}
{\rm div}\left(\si\nabla u\right)=0\quad\mbox{in }\RR^d.
\eeq
The gradient field $\nabla u$ is then said to be {\em isotropically realizable}.
{This reconstruction problem has been widely studied in the literature in terms of uniqueness, stability or instability, and algorithms of approximate solution (see, {\em e.g.}, \cite{FaEl}, \cite{Kno} and the references therein).}
The isotropy constraint is actually appropriate in Materials Science, since composite materials are built from isotropic phases. Moreover, the homogeneous conductivity equation \eqref{conequ} is satisfied by the local electric fields in periodic composites. We have proved in \cite{BMT} that any gradient field $\nabla u$ which is non-vanishing and regular is isotropically realizable in $\RR^d$. The main ingredient of this construction is the associated gradient flow
\beq\label{Xtxi}
\left\{\ba{rl}
\dis {\partial X\over \partial t}(t,x) & =\nabla u\big(X(t,x)\big)
\\ \ecart
X(0,x) & =x.
\ea\right.
\quad\mbox{for }t\in\RR,\ x\in\RR^d.
\eeq
The dynamical approach of \cite{BMT} forces the regularity {$u\in C^3(\RR^d)$}. However, this smoothness is not compatible with most of composite materials where the gradient is only piecewise regular (for instance regular in each phase of the material).
The purpose of the present work is to extend the results of \cite{BMT} to less regular gradient fields. To this end, we study two independent cases which are respectively developed in Section~\ref{s.nojump} and Section~\ref{s.jump}.
\par
In Section~\ref{s.nojump} we assume that the gradient field $\nabla u$ is continuous in $\RR^d$. The idea is to modify the strategy of \cite{BMT} applying the celebrated approach of DiPerna and Lions \cite{DiLi} for solving ordinary differential equations in suitable Sobolev spaces.
More precisely, we prove (see Theorem~\ref{thm.nojump} below) that any gradient field $\nabla u$ in $W^{1,1}_{\rm loc}(\RR^d)^d$ is isotropically realizable in $\RR^d$ if
\beq\label{Lapui}
\nabla u\mbox{ is uniformly continuous in }\RR^d,\quad \De u\in L^\infty(\RR^d)\quad\mbox{and}\quad\inf_{\RR^d}|\nabla u|>0.
\eeq
Moreover, any positive function $\si\in L^\infty_{\rm loc}(\RR^d)$ with $\si^{-1}\in L^\infty_{\rm loc}(\RR^d)$ is shown to be a suitable conductivity if and only if roughly speaking (see Remark~\ref{rem.continuity}) there exists $E$, a set of Lebesgue measure zero, such that
\beq\label{siXi}
{\sigma(x)\over\sigma\big(X(t,x)\big)}=\exp\left(\int_0^t \Delta u\big(X(s,x)\big)\,ds\right),\quad\forall\,t\in\RR,\ \forall\,x\in\RR^d\setminus E,
\eeq
where $X(\cdot,x)$ is the gradient flow \eqref{Xtxi}.
Assumption~\eqref{Lapui} improves significantly the regularity $u\in C^3(\RR^d)$ which is needed in~\cite{BMT}.
But the price to pay is that the reconstruction of an appropriate conductivity is much more delicate. In particular, by \cite{DiLi} the flow $X(\cdot,x)$ of~\eqref{Xtxi} is only continuous for almost everywhere $x\in\RR^d$. However, condition~\eqref{Lapui} is not still satisfactory since it excludes most of the Lipschitz continuous potentials~$u$ which naturally arise in composite materials.
\par
{In section~\ref{s.jump} we study the case of a piecewise regular gradient $\nabla u$ in a domain $\Om$ of $\RR^d$ composed by $n$ ``generalized" polyhedra $\Om_k$ ({\em i.e.} obtained from polyhedra through a smooth diffeomorphism). The continuous potential $u$ agrees in each set $\Om_k$ to a function $u_k\in C^2(\overline{\Om_k})$ such that the trajectories of \eqref {Xtxi} flow from an {\em inflow} boundary face (on which the outer normal derivative of $u_k$ is negative) to an {\em outflow} boundary face (on which the outer normal derivative of $u_k$ is positive), while the other boundary faces are tangential to $\nabla u_k$ (see Figure~\ref{fig1}). We prove (see Theorem~\ref{thm.jump} below) that there exists a piecewise continuous conductivity $\si$ solution to equation \eqref{conequ} if and only if for any contiguous polyhedra $\Om_j$ and $\Om_k$ of $\Om$, the normal derivatives satisfy the condition
\beq\label{ujuk}
{\partial u_j\over\partial\nu}={\partial u_k\over\partial\nu}=0\;\;\mbox{on }\partial\Om_j\cap\partial\Om_k\quad\mbox{or}\quad
{\partial u_j\over\partial\nu}\,{\partial u_k\over\partial\nu}>0\;\;\mbox{on }\partial\Om_j\cap\partial\Om_k.
\eeq
In the first case the common boundary face $\partial\Om_j\cap\partial\Om_k$ is tangential to the gradient, while in the second case $\partial\Om_j\cap\partial\Om_k$ is an inflow (resp. outflow) face of $\Om_j$ and an outflow (resp. inflow) face of $\Om_k$.
Actually, the picture is a little more constrained: We need to consider a so-called {\em $\nabla u$-admissible} domain~$\Om$ (see Definition~\ref{def.adm} below). Figure~\ref{fig2} below represents a $\nabla u$-admissible set, and Figure~\ref{fig3} represents a non-admissible one.
\par
We construct step by step a suitable piecewise conductivity $\si$ such that $\si=\si_k$ in $\Om_k$ as follows. If $\si_j$ is already constructed in $\Om_j$, by \cite{BoVa} and \cite{Ric} (see Proposition~\ref{pro.BVR} for details) there exists a unique positive function $\si_k\in C^1(\overline{\Om_k})$ solution to the equation ${\rm div}\left(\si_k\nabla u_k\right)=0$ in $\Om_k$, and equal on the inflow or outflow face $\partial\Om_j\cap\partial\Om_k$ to the boundary value $\ga_k\in C\big(\partial\Om_j\cap\partial\Om_k\big)$ which ensures by virtue of \eqref{ujuk} the flux continuity condition
\beq\label{conflu}
\si_j\,{\partial u_j\over\partial\nu}=\ga_k\,{\partial u_k\over\partial\nu}\quad\mbox{on }\partial\Om_j\cap\partial\Om_k.
\eeq
So, the piecewise continuous function $\si=\si_k$ in $\Om_k$ is a solution to the equation ${\rm div}\left(\si\nabla u\right)=0$ in the distributional sense of $\Om$.
}
\par
In Section~\ref{s.exa} the results of Section~\ref{s.jump} are illustrated by the case of piecewise constant gradients in some triangulation (see Figure~\ref{fig4} below), and the case of the gradient of a function $u\in C(\RR^d)$ defined by $u(x):=g_\pm(x_1)+f(x_2,\dots,x_d)$ in each half-space $\{\pm\,x_1>0\}$.
\subsubsection*{Notation}
\begin{itemize}
\item ${\rm int}\left(A\right)$ denotes the interior of a subset $A$ of $\RR^d$.
\item $C(A)$ denotes the set of continuous functions in a topological space $A$.
\item $C^k(A)$ denotes the space of $k$-differentiable functions in a subset $A$ of $\RR^d$, and $C^k_c(A)$ denotes the subspace of $C^k(A)$ composed of functions with compact support in $A$.
\item $\D'(\Om)$ denotes the distributions space in an open set $\Om$ of $\RR^d$.
\item $c$ denotes a positive constant which may vary from line to line.
\end{itemize}
\section{Case where the gradient field is continuous}\label{s.nojump}
For $u\in W^{2,1}_{\rm loc}(\RR^d)$, the gradient flow $X=X(t,x)$ associated with $\nabla u$ is defined (if possible) by
\beq\label{Xtx}
\left\{\ba{rl}
\dis {\partial X\over \partial t}(t,x) & =\nabla u\big(X(t,x)\big)
\\ \ecart
X(0,x) & =x.
\ea\right.
\quad\mbox{for }t\in\RR,\ x\in\RR^d.
\eeq
\begin{Thm}\label{thm.nojump}
Let $u:\RR^d\to\RR$ be a function satisfying
\beq\label{nojumpDu}
u\in W^{2,1}_{\rm loc}(\RR^d),\quad \nabla u\mbox{ is uniformly continuous in }\RR^d,\quad \Delta u\in L^\infty(\RR^d),
\quad \inf_{\RR^d}|\nabla u|>0.
\eeq
Then, there exists a positive function $\sigma\in L^\infty_{\rm loc}(\RR^d)$ with $\sigma^{-1}\in L^\infty_{\rm loc}(\RR^d)$, solution to the conductivity equation
\beq\label{divsiDu=0}
\div\left(\si\nabla u\right)=0\quad\mbox{in }\D'(\RR^d),
\eeq
the flow $X(\cdot,x)$ is well defined by \eqref{Xtx} for a.e. $x\in\RR^d$, and $\sigma$ satisfies
{the following: for any $t\in\RR$, there exists a set $E_t$, of Lebesgue measure zero depending on $t$, such that
\beq\label{siX}
{\sigma(x)\over\sigma\big(X(t,x)\big)}=\exp\left(\int_0^t \Delta u\big(X(s,x)\big)\,ds\right),\quad\forall\, x\in\RR^d\setminus E_t.
\eeq
Conversely, if there exists $E$, a set of Lebesgue measure zero, and a positive function $\sigma$ in $L^\infty_{\rm loc}(\RR^d)$ such that
\beq\label{siX2}
{\sigma(x)\over\sigma\big(X(t,x)\big)}=\exp\left(\int_0^t \Delta u\big(X(s,x)\big)\,ds\right)
\eeq
holds for any $t\in\RR$ and any $x\in\RR^d\setminus E$, then $\sigma$ is solution to equation~\eqref{divsiDu=0}.
}
\end{Thm}
\begin{Rem}
Assumptions \eqref{nojumpDu} replace the smoothness $u\in C^3(\RR^d)$ which is needed in~\cite{BMT}.
\end{Rem}
\begin{Rem}\label{rem.continuity}
{The set $E$ of Lebesgue measure zero where formula \eqref{siX2} is not satisfied by $x$ does not depend on $t$, while the set $E_t$ does depend on $t$ in formula \eqref{siX}. Hence, formula \eqref{siX2} is stronger than \eqref{siX}. Both formulas are equivalent if for instance $X$, $\De u$ and $\sigma$ are continuous.}
\end{Rem}
\noindent
{\bf Proof of Theorem~\ref{thm.nojump}.}
Let $(\rho_n)_{n\geq 1}$ be a sequence of mollifiers satisfying
\beq\label{rhon}
\rho_n\in C^\infty(\RR^d),\quad {\rm supp}\left(\rho_n\right)\subset B(0,1/n),\quad \rho_n\geq 0,\quad \int_{\RR^d}\rho_n(x)\,dx=1.
\eeq
Denote $u_n:=\rho_n*u\in C^\infty(\RR^d)$. Since by \eqref{nojumpDu} $\nabla u$ is uniformly continuous in $\RR^d$, the sequence $\nabla u_n=\rho_n*\nabla u$ converges uniformly to $\nabla u$ in $\RR^d$. Hence, by the last inequality of \eqref{nojumpDu} there exists a constant $m>0$ such that
\beq\label{Dun­0}
\inf_{\RR^d}|\nabla u_n|\geq m>0\quad\mbox{for $n$ large enough}.
\eeq
Let $X_n(t,x)$ be the flow associated with $\nabla u_n$ defined by
\beq\label{Xntx}
\left\{\ba{rl}
\dis {\partial X_n\over \partial t}(t,x) & =\nabla u_n\big(X(t,x)\big)
\\ \ecart
X_n(0,x) & =x.
\ea\right.
\quad\mbox{for }t\in\RR,\ x\in\RR^d.
\eeq
By \eqref{Dun­0} the regular case of \cite[Theorem~2.15]{BMT} shows that there exists a unique function $\tau_n$ in $C^\infty(\RR^d)$ satisfying
\beq\label{taun}
u_n\big(X_n(\tau_n(x),x)\big)=0,\quad\forall\,x\in\RR^d,
\eeq
and that, denoting
\beq\label{sigman}
\si_n(x):=\exp\left(\int_0^{\tau_n(x)}\Delta u_n\big(X_n(s,x)\big)\,ds\right)\quad\mbox{for }x\in\RR^d,
\eeq
we have
\beq\label{divsinDun=0}
\div\left(\si_n\nabla u_n\right)=0\quad\mbox{in }\RR^d,
\eeq
and
\beq\label{sinXn}
{\sigma_n(x)\over\sigma_n\big(X_n(t,x)\big)}=\exp\left(\int_0^t \Delta u_n\big(X_n(s,x)\big)\,ds\right),\quad\forall\,x\in\RR^d,\ \forall\,t\in\RR.
\eeq
\par
The main difficulty is now to pass to the limit $n\to\infty$ in equations \eqref{sigman}, \eqref{divsinDun=0}, \eqref{sinXn}.
To this end, we will use the approach of DiPerna and Lions \cite{DiLi} for solving ordinary differential equations in Sobolev spaces.
{First of all, note that by condition \eqref{nojumpDu} the field $b:=\nabla u$ satisfies the condition (49) and (70) of \cite{DiLi}, {\em i.e.}
\beq\label{conDu}
{b\over 1+|x|}\in L^\infty(\RR^d),\quad b\in W^{1,1}_{\rm loc}(\RR^d)^d\quad\mbox{and}\quad{\rm div}\,b\in L^\infty(\RR^d),
\eeq
since any uniformly continuous function $f(x)$ in $\RR^d$ is bounded by an affine function of $|x|$.
Hence, by virtue of \cite[Theorem III.2]{DiLi},} the flow $X_n(\cdot,x)$ converges in $C_{\rm loc}(\RR)$ to the unique flow {$X(\cdot,x)\in C^1(\RR^d)^d$ defined by \eqref{Xtx} for a.e. $x\in\RR^d$}. Moreover, $X$ satisfies the semi-group property: {for any $t\in\RR$, there exists a set $E_t$, of Lebesgue measure zero depending on $t$, such that
\beq\label{sgX}
X(s+t,x)=X\big(s,X(t,x)\big),\quad\forall\, s\in\RR,\ \forall\,x\in\RR^d\setminus E_t.
\eeq
}
{The image measure $\la_X(t)$, for $t\in\RR$, of the Lebesgue measure $\la$ by $X(t,\cdot)$, {\em i.e.} defined by
\beq\label{laX}
\int_{\RR^d}\ph\,d\la_X(t)=\int_{\RR^d}\ph\big(X(t,x)\big)\,dx,\quad\forall\,\ph\in C_c(\RR^d),
\eeq
has a density in $r(t,\cdot)\in L^\infty(\RR^d)$ with respect to the Lebesgue measure, which satisfies for any $t\in\RR$,
\beq\label{rho}
e^{-t\,\|\De u\|_{L^\infty(\RR^d)}}\leq r(t,\cdot)\leq e^{t\,\|\De u\|_{L^\infty(\RR^d)}}\;\;\mbox{a.e. in }\RR^d,
\eeq
or equivalently, for any $t\in\RR$ and for any $\ph\in C_c(\RR^d)$, $\ph\geq 0$,
\beq\label{elaX}
e^{-t\,\|\De u\|_{L^\infty(\RR^d)}}\int_{\RR^d}\ph(x)\,dx\leq\int_{\RR^d}\ph\,d\la_X(t)\leq e^{t\,\|\De u\|_{L^\infty(\RR^d)}}\int_{\RR^d}\ph(x)\,dx.
\eeq
}
\par
We will need the following result satisfied by the flows $X_n$ and $X$.
\begin{Lem}\label{lem.flow}
\hfill\noindent
\begin{itemize}
{\item[$i)$] If $f\in L^1_{\rm loc}(\RR^d)$ then $f\circ X\in L^1_{\rm loc}(\RR\times\RR^d)$.}
\item[$ii)$] Let $f\in L^1_{\rm loc}(\RR^d)$, let $K$ be a compact of $\RR^d$, and let $I$ be a bounded interval of $\RR$. Then, we have
\beq\label{fXnX}
\lim_{n\to\infty}\int_K\int_I\left|f\big(X_n(s,x)\big)-f\big(X(s,x)\big)\right| ds\,dx=0.
\eeq
\item[$iii)$] Let $f_n$ be a non-negative sequence of $L^1_{\rm loc}(\RR^d)$ which converges strongly to $0$ in $L^1_{\rm loc}(\RR^d)$,
let $K$ be a compact of $\RR^d$, and let $I$ be a bounded interval of $\RR$. Then, we have
\beq\label{fnXn}
\lim_{n\to\infty}\int_K\int_I f_n\big(X_n(s,x)\big)\,ds\,dx=0.
\eeq
\item[$iv)$] Let $F\in L^p(\RR^d)^N$ for $N\in\NN$, $p\in[1,\infty)$, let $G\in L^{p'}(\RR^d)^N$ with compact support, where $p'$ is the conjugate exponent of $p$, and let $\rho_n$ be a sequence in $C^\infty_c(\RR)$ satisfying \eqref{rhon} with $d=1$. Then, we have
\beq\label{rhonFXG}
\lim_{n\to\infty}\int_{\RR^d}\int_\RR \rho_n(s)\,F\big(X(s,x)\big)\cdot G(x)\,ds\,dx=\int_{\RR^d} F(x)\cdot G(x)\,dx.
\eeq
\end{itemize}
\end{Lem}
The proof is divided in {five} steps.
\par\ms\noindent
{\it First step:} Convergence of the sequence $\tau_n$ defined by \eqref{taun}.
\par\noindent
On the one hand, since by \eqref{nojumpDu} {there exists $E$, a set of Lebesgue measure zero, such that for any $x\in\RR\setminus E$,
\[
{d\over dt}\big(u(X(t,x)\big)=|\nabla u|^2\big(X(t,x)\big)\geq\inf_{\RR^d}|\nabla u|^2>0,\quad\forall\,t\in\RR,
\]
}
there exists a unique $\tau(x)\in\RR$ such that
\beq\label{tau}
u\big(X(\tau(x),x)\big)=0\quad\mbox{for a.e. }x\in\RR^d.
\eeq
On the other hand, by \eqref{taun} we have
\beq\label{untaun}
|u_n(x)|=\left|u_n(x)-u_n\big(X_n(\tau_n(x),x)\big)\right|=\left|\,\int_0^{\tau_n(x)}|\nabla u_n|^2\big(X_n(t,x)\big)\,ds\,\right|\geq m^2\,|\tau_n(x)|
\quad\mbox{a.e. }x\in\RR^d.
\eeq
Hence, since $u_n$ converges uniformly to $u$ in any compact set $K$ of $\RR^d$, the sequence $\tau_n$ is bounded in $L^\infty(K)$.
Let $x\in\RR^d$ be satisfying \eqref{untaun}. Up to a subsequence still denoted by~$n$, $\tau_n(x)$ converges to some $\tau_x$ in $\RR$.
Using the uniform convergence of $X_n(\cdot,x)$ to $X(\cdot,x)$ and passing to the limit in equality \eqref{taun} we get that $u\big(X(\tau_x,x)\big)=0$, which by uniqueness of $\tau(x)$ implies that $\tau_x=\tau(x)$.
Therefore, we obtain for the whole sequence
\beq\label{contaun}
\lim_{n\to\infty}\tau_n(x)=\tau(x)\quad\mbox{for a.e. }x\in\RR^d.
\eeq
\par
{Since $\tau$ is measurable and $\De u\circ X\in L^1_{\rm loc}(\RR\times\RR^d)$ by Lemma~\ref{lem.flow}, applying Fubini's theorem to the function $(t,x)\mapsto 1_{[0,\tau(x)]}(t)\,\De u\big(X(t,x)\big)$ in $L^1_{\rm loc}(\RR\times\RR^d)$, we can define the measurable function $\si$ by
}
\beq\label{sitau}
\si(x):=\exp\left(\int_0^{\tau(x)}\De u\big(X(s,x)\big)\,ds\right)\quad\mbox{for a.e. }x\in\RR^d. 
\eeq
{\it Second step:} Strong convergence of the sequence $w_n:=\ln\sigma_n$ to $w:=\ln\sigma$ in $L^1_{\rm loc}(\RR^d)$.
\par\smallskip\noindent
Let $K$ be a compact set of $\RR^d$. We have
\beq\label{Eni}
\ba{rll}
\dis \int_K |w_n(x)-w(x)|\,dx\leq
& \dis \int_K \left|\,\int_0^{\tau_n(x)}\left|\De u\big(X_n(s,x)\big)-\De u\big(X(s,x)\big)\right|ds\,\right|dx & =:E_n^1
\\ \ecart
+ & \dis \int_K\left|\,\int_0^{\tau_n(x)}\big|\De u_n-\De u\big|\big(X_n(s,x)\big)\,ds\,\right|dx & =:E_n^2
\\ \ecart
+ & \dis \int_K \left|\,\int_{\tau(x)}^{\tau_n(x)}\left|\De u\big(X(s,x)\big)\right|ds\,\right|dx & =:E_n^3.
\ea
\eeq
Since by the first step the sequence $\tau_n$ is uniformly bounded in any compact set of $\RR^d$, there exist a bounded interval $I$ of $\RR$ such that
\[
E_n^1\leq \int_K\int_I\left|\De u\big(X_n(s,x)\big)-\De u\big(X(s,x)\big)\right|ds\,dx.
\]
Hence, applying the limit \eqref{fXnX} of Lemma~\ref{lem.flow} with $f:=\De u$, we get that $E_n^1$ tends to $0$.
Similarly, applying \eqref{fnXn} with the sequence $f_n:=\De u_n-\De u=\rho_n*\De u-\De u$ which converges strongly to $0$ in $L^1_{\rm loc}(\RR^d)$, we get that $E_n^2$ tends to $0$. Finally, since $\tau_n$ is uniformly bounded in the compact $K$ and $\De u\in L^\infty(\RR^d)$, by convergence \eqref{contaun} and the Lebesgue dominated convergence theorem we get that
\[
0\leq E_n^3\leq c\int_K|\tau_n-\tau|\,dx\limi 0.
\]
Therefore, passing to the limit $n\to\infty$ in \eqref{Eni} we obtain that the sequence $w_n$ converges strongly to $w$ in $L^1_{\rm loc}(\RR^d)$.
\par\smallskip\noindent
{\it Third step:} Derivation of the conductivity equation \eqref{divsiDu=0}.
\par\smallskip\noindent
By \eqref{sigman} the function $w_n$ is defined by
\beq\label{wn}
w_n(x)=\int_0^{\tau_n(x)}\Delta u_n\big(X_n(s,x)\big)\,ds\quad\mbox{for }x\in\RR^d.
\eeq
Since by the first step $\tau_n$ is bounded in any compact of $\RR^d$ and $\De u_n=\rho_n*\De u$ is bounded in $L^\infty(\RR^d)$, the sequence $w_n$ is bounded in $L^\infty_{\rm loc}(\RR^d)$. Hence, by the second step the sequence $\sigma_n=e^{w_n}$ converge strongly to $\sigma=e^w$ in $L^1_{\rm loc}(\RR^d)$. Moreover, the sequence $\nabla u_n$ converges to $\nabla u$ in $C_{\rm loc}(\RR^d)$. Therefore, passing to the limit in equation \eqref{divsinDun=0} we get that $\si$ is solution to the conductivity equation \eqref{divsiDu=0} in the distributions sense. Finally, both $\si$ and $\si^{-1}$ belong to $L^\infty_{\rm loc}(\RR^d)$, since $\si$ is the limit in $L^1_{\rm loc}(\RR^d)$ of the sequence $\si_n=e^{w_n}$ which is bounded in $L^\infty_{\rm loc}(\RR^d)$.
\par\smallskip\noindent
{\it Fourth step:} Proof of formula \eqref{siX}.
\par\smallskip\noindent
Formula \eqref{sinXn} reads as
\beq\label{fwn}
w_n(x)-w_n\big(X_n(t,x)\big)=\int_0^t\De u_n\big(X_n(s,x)\big)\,ds,\quad\forall\,t\in\RR,\ \forall\,x\in\RR^d.
\eeq
On the one hand, writing
\[
\left|w_n\big(X_n(t,x)\big)-w\big(X(t,x)\big)\right|\leq \left|w\big(X_n(t,x)\big)-w\big(X(t,x)\big)\right|+|w_n-w|\big(X_n(t,x)\big),
\]
applying limit \eqref{fXnX} with $f:=w$, and applying limit \eqref{fnXn} with $f_n:=|w_n-w|$ which converges strongly to $0$ in $L^1_{\rm loc}(\RR^d)$ by the second step, we get that
\beq\label{wnXnwX}
w_n\big(X_n(t,\cdot)\big)\limi w\big(X(t,\cdot)\big)\quad\mbox{strongly in }L^1_{\rm loc}(\RR^d),\ \mbox{for any }t\in\RR.
\eeq
On the other hand, let $K$ be a compact set of $\RR^d$ and $t\in\RR$. We have
\[
\ba{l}
\dis \int_K\left|\,\int_0^t\De u_n\big(X_n(s,x)\big)\,ds-\int_0^t\De u\big(X(s,x)\big)\,ds\,\right|dx
\\ \ecart
\dis \leq \left|\,\int_0^t\int_K\Big[\left|\De u\big(X_n(s,x)\big)-\De u\big(X(s,x)\big)\right|
+\big|\De u_n-\De u\big|\big(X_n(s,x)\big)\Big]dx\,ds\,\right|.
\ea
\]
{Then, applying successively limit \eqref{fXnX} with $f:=\De u$ and limit \eqref{fnXn} with $f_n:=|\De u_n-\De u|$ in $[0,t]\times K$, we get that}
\beq\label{LunXnLuX}
\int_0^t\De u_n\big(X_n(s,x)\big)\,ds\limi\int_0^t\De u\big(X(s,x)\big)\,ds
\quad\mbox{strongly in }L^1_{\rm loc}(\RR^d),\ \mbox{for any }t\in\RR.
\eeq
Therefore, using the limits \eqref{wnXnwX} and \eqref{LunXnLuX} in \eqref{fwn}, {there exists $E_t$, a set of Lebesgue measure zero depending on $t$, such that for any $t\in\RR$,
\beq\label{fw}
w(x)-w\big(X(t,x)\big)=\int_0^t\De u\big(X(s,x)\big)\,ds,\quad\forall\,x\in\RR^d\setminus E_t.
\eeq
}
or equivalently formula \eqref{siX}.
\begin{Rem}\label{rem.proofsi}
A direct proof of \eqref{siX} would consist in replacing $x$ by $X(t,x)$ in the definition \eqref{sitau} of $\si(x)$ and to use the semi-group property \eqref{sgX}, to obtain the desired formula \eqref{siX}. However, since the function $\tau$ involving in \eqref{sitau} is only defined a.e. in $\RR^d$ by \eqref{tau}, it is not clear that for an admissible point~$x$ of $\tau$, $X(t,x)$ for $t\in\RR$, is also an admissible point of $\tau$.
\end{Rem}
\par\smallskip\noindent
{\it Fifth step:} Formula \eqref{siX2} implies the conductivity equation \eqref{divsiDu=0}.
\par\smallskip\noindent
Let $\si$ be a positive function in $L^\infty_{\rm loc}(\RR^d)$ satisfying formula \eqref{siX}.
First of all by \eqref{nojumpDu} the function $b(t,x):=\nabla u(x)$ satisfies the assumptions $(*),(**)$ of \cite[Theorem~II.3.1)]{DiLi} and assumptions $(49),(70)$ of \cite[Theorem~III.2]{DiLi}. Then, by virtue of \cite[Theorem~II.3.1)]{DiLi} and \cite[Theorem~III.2]{DiLi} the function $\sigma\big(X(t,x)\big)$ is solution to the transport equation
\beq\label{transiX}
{\partial\over\partial t}\left[\sigma\big(X(t,x)\big)\right]=\nabla u(x)\cdot\nabla_x\left[\sigma\big(X(t,x)\big)\right]\quad\mbox{in }\D'(\RR\times\RR^d).
\eeq
Moreover, taking the derivative with respect to $t$ in \eqref{siX2} (at this point \eqref{siX} seems to be not sufficient) we have
\[
{\partial\over\partial t}\left[\sigma\big(X(t,x)\big)\right]=-\,\sigma\big(X(t,x)\big)\,\De u\big(X(t,x)\big)\quad\mbox{in }\D'(\RR\times\RR^d).
\]
Equating the two previous equations we get that
\[
\nabla_x\left[\sigma\big(X(t,x)\big)\right]\cdot\nabla u(x)+\sigma\big(X(t,x)\big)\,\De u\big(X(t,x)\big)=0\quad\mbox{in }\D'(\RR\times\RR^d).
\]
Since $\nabla u\in W^{1,1}_{\rm loc}(\RR^d)$, the previous equation can be read as
\[
{\rm div}_x\left[\sigma\big(X(t,x)\big)\nabla u(x)\right]=\sigma\big(X(t,x)\big)\left[\De u(x)-\De u\big(X(t,x)\big)\right]\quad\mbox{in }\D'(\RR\times\RR^d),
\]
which implies that for any $\ph\in C^\infty_c(\RR)$ and $\psi\in C^\infty_c(\RR^d)$,
\beq\label{divsiXDu}
\ba{l}
\dis \int_{\RR^d}\int_\RR \ph(t)\,\sigma\big(X(t,x)\big)\nabla u(x)\cdot\nabla\psi(x)\,dt\,dx
\\ \ecart
\dis =\int_{\RR^d}\int_\RR \ph(t)\,\psi(x)\,\sigma\big(X(t,x)\big)\left[\De u\big(X(t,x)\big)-\De u(x)\right]dt\,dx.
\ea
\eeq
Taking $\ph(t)=\rho_n(t)$ in \eqref{divsiXDu} and applying the limit \eqref{rhonFXG} of Lemma~\ref{lem.flow} with $F=\si,\si,\si\De u$ in $L^{p}_{\rm loc}(\RR^d)$ for $p:={d\over d-1}$, and respectively $G=\nabla u\cdot\nabla\psi,\psi\De u,\psi$ in $L^{p'}(\RR^d)$ with compact support, we obtain that
\[
\int_{\RR^d}\sigma(x)\nabla u(x)\cdot\nabla\psi(x)\,dx=0,\quad\forall\,\psi\in C^\infty_c(\RR^d),
\]
or equivalently the conductivity equation \eqref{divsiDu=0}.
\par\bs\noindent
{\bf Proof of Lemma~\ref{lem.flow}.}
\par\smallskip\noindent
{
$i)$ Let $I$ be a bounded interval of $\RR$ and let $K$ be a compact set of $\RR^d$.
We have for any $t\in I$ and $x\in K$,
\[
\big|X_n(t,x)\big|\leq |x|+\left|\,\int_0^t \left|\nabla u_n\big(X_n(s,x)\big)\right|\,ds\,\right|.
\]
Moreover, the uniform continuity of $\nabla u$ in $\RR^d$ and the equality $\nabla u_n=\rho_n*\nabla u$ imply the existence of a constant $c>0$ such that
\[
|\nabla u_n(y)|\leq c\,|y|+c,\quad\forall\,n\in\NN,\ \forall\,y\in\RR^d.
\]
We thus deduce that
\[
\big|X_n(t,x)\big|\leq c+c\,\left|\,\int_0^t \left|X_n(s,x)\right|\,ds\,\right|,\quad\forall\,n\in\NN,\ \forall\,t\in I,\ \forall\,x\in K.
\]
Hence, by Gronwall's inequality (see, {\em e.g.}, \cite[Section 17.3]{HSD}) there exists a constant $c>0$ such that
\beq\label{estXn}
\big|X_n(t,x)\big|\leq c\,e^{c\,|t|},\quad\forall\,n\in\NN,\ \forall\,t\in I,\ \forall\,x\in K.
\eeq
Therefore, there exists a compact $\hat{K}$ of $\RR^d$ and $E$, a set of Lebesgue measure zero, such that
\beq\label{hK}
X_n(t,x),X(t,x)\in \hat{K},\quad\forall\,n\in\NN,\ \,\forall\,t\in I,\ \forall\,x\in K\setminus E.
\eeq
\par
Let $f\in L^1_{\rm loc}(\RR^d)$, and let $f_n$ be a sequence in $C^\infty_c(\RR^d)$ which converges strongly to $f$ in $L^1_{\rm loc}(\RR^d)$.
We will show that $f_n\circ X$ converges strongly to some function $g$ in $L^1(I\times K)$.
By \cite[Theorem~II.3.1)]{DiLi} and \cite[Theorem~III.2]{DiLi} $f_n\circ X$ is in $L^1_{\rm loc}(\RR\times\RR^d)$.
Let $O$ be a bounded open set of $\RR^d$ containing the compact set $\hat{K}$, and let $\psi$ be a non-negative function in $C_c(O)$ which is equal to $1$ in $\hat{K}$.
By \eqref{hK} and estimate \eqref{elaX} we have for any $p,q\in\NN$,
\[
\ba{ll}
\dis \int_I\int_K\big|f_p(X(t,x))-f_q(X(t,x))\big|\,dt\,dx & \dis \leq\int_I dt\int_{\RR^d}\psi(X(t,x))\,\big|f_p(X(t,x))-f_q(X(t,x))\big|\,dx
\\ \ecart
& \dis =\int_I dt\int_{\RR^d}\psi\,|f_p-f_q|\,d\la_{X}(t)\leq c\int_{O}|f_p-f_q|.
\ea
\]
Hence, $f_n\circ X$ is a Cauchy sequence in $L^1(I\times K)$ and thus converges strongly to some function $g$ in $L^1(I\times K)$. Therefore, due to the arbitrariness of $I,K$ the sequence $f_n\circ X$ converges strongly to some function $g$ in $L^1_{\rm loc}(\RR\times\RR^d)$.
\par
Finally, by estimate \eqref{rho} we have for any bounded interval $I$ of $\RR$, any bounded open set $O$ of $\RR^d$ and any function $\ph\in C_c(O)$,
\[
\ba{l}
\dis \int_I dt\int_{\RR^d}\ph f \,d\la_X(t)=\int_I dt\int_{O}\ph(x)\,f(x)\,r(t,x)\,dx=\lim_{n\to\infty}\int_I dt\int_{O}\ph(x)\,f_n(x) \,r(t,x)\,dx
\\ \ecart
\dis =\lim_{n\to\infty}\int_I dt\int_{\RR^d}\ph f_n\,d\la_X(t)=\lim_{n\to\infty}\int_I \int_{\RR^d}(\ph f_n)\big(X(t,x)\big)\,dx=\int_I dt\int_{\RR^d}\ph\big(X(t,x)\big)\,g(x)\,dx,
\ea
\]
which, due to the arbitrariness of $I,O,\ph$, implies that $f\circ X=g\in L^1_{\rm loc}(\RR\times \RR^d)$.
}
\par\smallskip\noindent
$ii)$ Let $I$ be a bounded interval of $\RR$ and let $K$ be a compact set of $\RR^d$.
Let $\ph\in C^\infty_c(\RR^d)$ be an approximation of $f$ in $L^1(\RR^d)$.
We have
\beq\label{estfph1}
\ba{l}
\dis \limsup_{n\to\infty}\int_K\int_I\left|f\big(X_n(s,x)\big)-f\big(X(s,x)\big)\right| ds\,dx
\\ \ecart
\dis \leq\limsup_{n\to\infty}\int_K\int_I\left|\ph\big(X_n(s,x)\big)-\ph\big(X(s,x)\big)\right| ds\,dx
\\ \ecart
\dis +\limsup_{n\to\infty}\int_K\int_I |f-\ph|\big(X_n(s,x)\big)\,ds\,dx+\int_K\int_I |f-\ph|\big(X(s,x)\big)\,ds\,dx.
\ea
\eeq
On the one hand, the uniform convergence of $X_n(\cdot,x)$ to $X(\cdot,x)$ in $I$ combined with the continuity of $\ph$ yields that
\[
\int_I\left|\ph\big(X_n(s,x)\big)-\ph\big(X(s,x)\big)\right| ds\limi 0\quad\mbox{a.e. }x\in K,
\]
and estimate \eqref{estXn} combined with the continuity of $\ph$ gives that
\[
\int_I\left|\ph\big(X_n(s,x)\big)-\ph\big(X(s,x)\big)\right| ds\leq c\quad\mbox{a.e. }x\in K.
\]
Hence, by the Lebesgue dominated convergence theorem
\beq\label{estfph2}
\lim_{n\to\infty}\int_K\int_I\left|\ph\big(X_n(s,x)\big)-\ph\big(X(s,x)\big)\right| ds\,dx=0.
\eeq
Then, {since by \eqref{hK} there exists a set $E$, of Lebesgue measure zero, such that
\beq\label{1hK}
1_{K}(x)\leq\min\big(1_{\hat{K}}(X(t,x)),1_{\hat{K}}(X_n(t,x))\big),\quad \forall\,n\in\NN,\ \forall\,t\in I,\ \forall\,x\in\RR^d\setminus E,
\eeq
}
using the estimate \eqref{elaX} satisfied by the image measure $\la_X(s)$ with $\De u$ and the similar one satisfied by $\la_{X_n}(s)$ with $\De u_n$, we get that
\[
\ba{l}
\dis \limsup_{n\to\infty}\int_K\int_I |f-\ph|\big(X_n(s,x)\big)\,ds\,dx+\int_K\int_I |f-\ph|\big(X(s,x)\big)\,ds\,dx
\\ \ecart
\dis \leq \limsup_{n\to\infty}\int_I \int_{\RR^d}\big(1_{\hat{K}}|f-\ph|\big)\big(X_n(s,x)\big)\,dx\,ds
+\int_I \int_{\RR^d}\big(1_{\hat{K}}|f-\ph|\big)\big(X(s,x)\big)\,dx\,ds
\\ \ecart
\dis = \limsup_{n\to\infty}\int_I \int_{\RR^d}1_{\hat{K}}(y)\,|f-\ph|(y)\,\la_{X_n}(s)(dy)\,ds
+\int_I \int_{\RR^d}1_{\hat{K}}(y)\,|f-\ph|(y)\,\la_{X}(s)(dy)\,ds
\\ \ecart
\leq c\,\|f-\ph\|_{L^1(\hat{K})}.
\ea
\]
Therefore, putting this and limit \eqref{estfph2} in \eqref{estfph1} we deduce the desired limit \eqref{fXnX}.
\par\ms\noindent
$iii)$ Let $I$ be a bounded interval of $\RR$, let $K$ be a compact set of $\RR^d$, and let $\hat{K}$ be a compact set of $\RR^d$ satisfying \eqref{hK}. Let $f_n$ be a non-negative sequence of $L^1_{\rm loc}(\RR^d)$ which converges strongly to $0$ in $L^1_{\rm loc}(\RR^d)$. Repeating the argument of $ii)$ {using inequality \eqref{1hK}} and the estimate \eqref{elaX} with $X_n$ in place of~$X$, we get that 
\[
\ba{ll}
\dis \limsup_{n\to\infty}\int_K\int_I f_n\big(X_n(s,x)\big)\,ds\,dx
& \dis \leq\limsup_{n\to\infty}\int_I\int_{\RR^d}\big(1_{\hat{K}}f_n\big)\big(X_n(s,x)\big)\,ds\,dx
\\ \ecart
& \dis \leq\limsup_{n\to\infty}\int_I\int_{\RR^d}1_{\hat{K}}(y)\,f_n(y)\,\la_{X_n}(s)(dy)\,ds
\\ \ecart
& \dis \leq c\,\limsup_{n\to\infty}\|f_n\|_{L^1(\hat{K})}=0,
\ea
\]
which yields \eqref{fnXn}.
\par\ms\noindent
$iv)$ Let $F\in L^p_{\rm loc}(\RR^d)^N$ for $N\in\NN$, $p\in[1,\infty)$, and let $G\in L^{p'}(\RR^d)^N$ whose support is included in a compact set $K$ of $\RR^d$. Consider a compact set $\hat{K}$ of $\RR^d$ satisfying \eqref{hK} with $I=[-1,1]$ and $K$, {\em i.e.} {there exists a set $E$, of Lebesgue measure zero, such that
\[
1_{\hat{K}}\big(X(t,x)\big)=1,\quad\forall\,t\in[-1,1],\ \forall\,x\in K\setminus E.
\]
}
Let $\Phi\in C^\infty_c(\RR^d)^N$ be an approximation of $F$ in $L^p(\hat{K})^N$.
By \eqref{rhon} we have
\[
\ba{l}
\dis \int_{\RR^d}\int_\RR \rho_n(s)\,F\big(X(s,x)\big)\cdot G(x)\,ds\,dx-\int_{\RR^d} F(x)\cdot G(x)\,dx
\\ \ecart
\dis =\int_{\RR^d}\int_\RR \rho_n(s)\left[\Phi\big(X(s,x)-\Phi(x)\big)\right]\cdot G(x)\,ds
\\ \ecart
\dis +\int_{\RR^d}\int_\RR \rho_n(s)\left[\big(1_{\hat{K}}(F-\Phi)\big)\big(X(s,x)\big)-\big(1_{\hat{K}}(F-\Phi)\big)(x)\right]\cdot G(x)\,ds\,dx.
\\ \ecart
\ea
\]
Then, by the H\"older inequality combined with estimate \eqref{rho} we get that
\beq\label{estFG}
\ba{l}
\dis \limsup_{n\to\infty}\left|\,\int_{\RR^d}\int_\RR \rho_n(s)\,F\big(X(s,x)\big)\cdot G(x)\,ds\,dx-\int_{\RR^d} F(x)\cdot G(x)\,dx\,\right|
\\ \ecart
\dis \leq \limsup_{n\to\infty}\left|\,\int_{\RR^d}\int_\RR \rho_n(s)\left[\Phi\big(X(s,x)\big)-\Phi(x)\right]\cdot G(x)\,ds\,dx\,\right|
+c\,\|F-\Phi\|_{L^p(\hat{K})^N}\,\|G\|_{L^{p'}(K)^N}.
\ea
\eeq
By the continuity of $\Phi$ we have
\[
\int_\RR \rho_n(s)\left[\Phi\big(X(s,x)\big)-\Phi(x)\right]ds\limi 0\quad\mbox{a.e. }x\in\RR^d.
\]
Moreover, we have
\[
\left|\,\int_\RR \rho_n(s)\left[\Phi\big(X(s,x)\big)-\Phi(x)\right]ds\,\right|\leq 2\,\|\Phi\|_{L^\infty(\RR^d)^N}\quad\mbox{a.e. }x\in\RR^d,
\]
so that
\[
\left|\left(\int_\RR \rho_n(s)\left[\Phi\big(X(s,x)\big)-\Phi(x)\right]ds\,\right)\cdot G(x)\,\right|\leq c\,|G(x)|\quad\mbox{a.e. }x\in\RR^d.
\]
Hence, since $G\in L^1(\RR^d)^N$ due to its compact support, the Lebesgue dominated convergence theorem implies that
\[
\lim_{n\to\infty}\int_{\RR^d}\int_\RR \rho_n(s)\left[\Phi\big(X(s,x)\big)-\Phi(x)\right]\cdot G(x)\,ds\,dx=0.
\]
Using this in \eqref{estFG} we thus obtain limit \eqref{rhonFXG}.
\cqfd
{
\section{Case where the gradient field has jumps}\label{s.jump}
In this section we will consider a gradient field which is piecewise regular in a finite number of so-called gradient-admissible domains.
\subsection{Gradient-admissible domain}
The starting point is the following result first due to Bongiorno, Valente \cite{BoVa}, and well reformulated by Richter \cite{Ric}.
\begin{Pro}[{\cite[Lemma~2]{Ric}}]\label{pro.BVR}
Let $\Om$ be a bounded domain ({\em i.e.} a connected open set) of~$\RR^d$, and let $u\in C^2(\overline{\Om})$ such that
\beq\label{Du>0}
\inf_{\Om}|\nabla u|>0.
\eeq
Let $\Ga_-$ be the {\em inflow} boundary of $\Om$, {\em i.e.} the subset of $\partial\Om$ on which the outer normal derivative of $u$ is negative: ${\partial u\over\partial\nu}<0$, and let $\Ga_+$ be the {\em outflow} boundary of $\Om$, {\em i.e.} the subset of $\partial\Om$ on which the outer normal derivative of $u$ is positive: ${\partial u\over\partial\nu}>0$.
\par\noindent
Then, each point of $\Om$ belongs to a unique trajectory $t\mapsto X(t,x)$ which flows from $\Ga_-$ to $\Ga_+$.
Moreover, there exists a unique positive function $\si\in C^1(\overline{\Om})$ taking prescribed values on $\Ga_-$ (resp. on $\Ga_+)$  which is solution to the equation
${\rm div}\left(\si\nabla u\right)=0$ in $\Om$.
\end{Pro}
\begin{Rem}
Actually, in \cite{Ric} the existence and the uniqueness of the conductivity $\si$ taking previous values on the inflow boundary $\Ga_-$ is proved under the weaker assumption
\[
\inf_{\Om}\big(\min\,(|\nabla u|,\De u)\big)>0.
\]
However, we will need the stronger condition \eqref{Du>0} in the sequel.
\end{Rem}
\noindent
{\bf Proof of Proposition~\ref{pro.BVR}.}
The proof can be found in \cite{Ric}. We will give another expression of the conductivity $\si$ following Theorem~\ref{thm.nojump}.
Let $\ga$ be a positive function in $C^1(\overline{\Ga_-})$. For a fixed $x\in\Om$, the trajectory $t\in[\tau_-(x),\tau_+(x)]\mapsto X(t,x)$ flows from the inflow boundary $\Ga_-$ to the outflow boundary $\Ga_+$, where $\tau_-(x)<0<\tau_+(x)$ and $X\big(\tau_\pm(x),x\big)\in\Ga_\pm$.
Let $y=X(\tau,x)$ be a point on the same trajectory. Note that by the semi-group property of the flow we have
\[
X\big(\tau_-(x),x\big)=X\big(\tau_-(y),y\big)=X\big(\tau_-(y),X(\tau,x)\big)=X\big(\tau_-(y)+\tau,x\big),
\]
hence $\tau_-(y)=\tau_-(x)-\tau$. Now, we can define the conductivity $\si_\ga$ along the trajectory by
\beq\label{sigatau}
\sigma_\ga\big(X(t,x)\big):=\ga\big(X(\tau_-(x),x)\big)\,\exp\left(\int_t^{\tau_-(x)}\Delta u\big(X(s,x)\big)\,ds\right)
\quad\mbox{for }t\in[\tau_-(x),\tau_+(x)].
\eeq
Formula \eqref{sigatau} does not depend on the point $y=X(\tau,x)$ on the same trajectory, since
\[
\int_{\tau_-(y)}^t \Delta u\big(X(s,y)\big)\,ds=\int_{\tau_-(x)-\tau}^t \Delta u\big(X(s+\tau,x)\big)\,ds
=\int_{\tau_-(x)}^{t+\tau}\Delta u\big(X(s,x)\big)\,ds,
\]
which implies that $\si_\ga\big(X(t,y)\big)=\si_\ga\big(X(t+\tau,x)\big)$.
Moreover, it is immediate that formula \eqref{sigatau} implies formula \eqref{siX2}.
Therefore, by Theorem~\ref{thm.nojump} $\si_\ga$ is a solution to the equation ${\rm div}\left(\si_\ga\nabla u\right)=0$ in $\Om$, and $\si_\ga=\ga$ on $\Ga_-$.
\par
Conversely, consider a positive function $\si\in C^1(\overline{\Om})$ such that ${\rm div}\left(\si\nabla u\right)=0$ in $\Om$, and $\si=\ga$ on~$\Ga_-$.
From the equality $\nabla\si\cdot\nabla u+\si\,\De u=0$ in $\Om$, we deduce that for any $x\in\Om$,
\[
{d\over dt}\left[\ln\big(\si(X(t,x)\big)\right]=-\De u\big(X(t,x)\big),\quad\forall\,t\in[\tau_-(x),\tau_+(x)],
\]
then
\[
{\si(x)\over\sigma\big(X(t,x)\big)}=\exp\left(\int_0^t \Delta u\big(X(s,x)\big)\,ds\right),\quad\forall\,t\in[\tau_-(x),\tau_+(x)].
\]
This combined with \eqref{sigatau} implies  that for any $x\in\Om$,
\[
{\sigma(x)\over\ga\big(X(\tau_-(x),x)\big)}={\sigma(x)\over\si\big(X(\tau_-(x),x)\big)}=\exp\left(\int_0^{\tau_-(x)}\Delta u\big(X(s,x)\big)\,ds\right)
={\sigma_\ga(x)\over\ga\big(X(\tau_-(x),x)\big)}.
\]
Therefore, we obtain that $\si=\si_\ga$ in $\Om$, which shows the uniqueness of the conductivity $\si_\ga$.
\cqfd
\par\bs
We can now state the definition of a gradient-admissible set.
\begin{Def}\label{def.Du-adm}
Let $\Om$ be a bounded domain of $\RR^d$, and let $u\in C^2(\overline{\Om})$. The domain $\Om$ is said to be {\em $\nabla u$-admissible} if condition \eqref{Du>0} holds.
\end{Def}
\begin{Rem}
The boundary of a $\nabla u$-admissible domain $\Om$ is split into the inflow boundary $\Ga_-$, the outflow boundary $\Ga_+$, and  surfaces which are tangential to $\nabla u$. Figure~\ref{fig1} shows a two-dimensional $\nabla u$-admissible domain $\Om$ with two boundary curves which are tangential to~$\nabla u$.
\end{Rem}
\begin{figure}
\centering
\vskip -1.cm
\includegraphics[scale=.6]{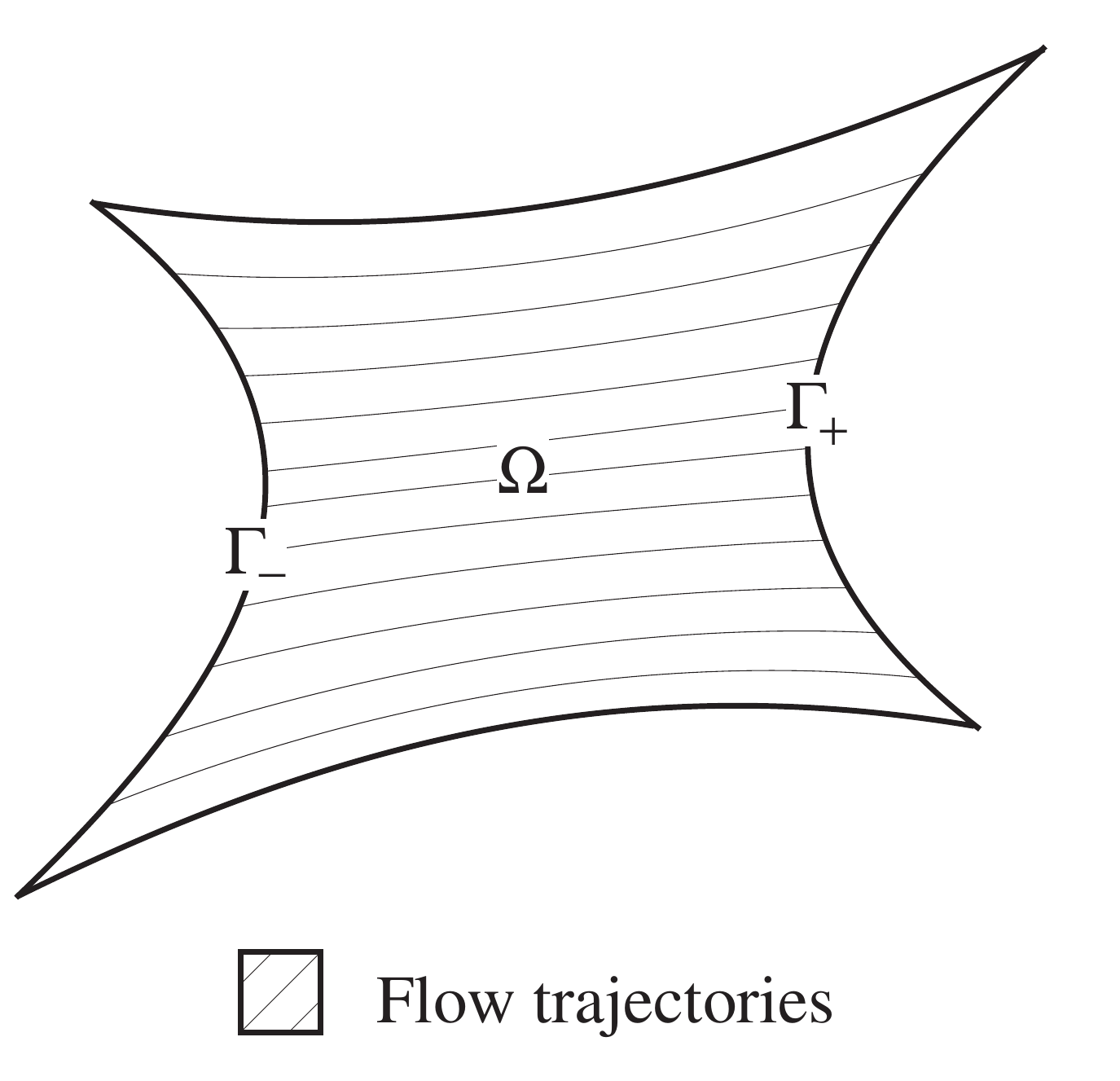}
\vskip 0.cm
\caption{\it The trajectories in $\Om$ flow from $\Gamma_-$ to $\Gamma_+$}
\label{fig1}
\end{figure}
\subsection{Piecewise regular gradient field}
In connection with the definition~\ref{def.Du-adm} of a gradient-admissible set, we focus on a so-called {\em admissible} domain defined as follows.
\begin{Def}\label{def.adm}
Let $\Om$ be a bounded domain of $\RR^d$. The set $\Om$ is said to be {\em admissible} if it is decomposed into  ``generalized open polyhedra" (obtained from polyhedra through a smooth diffeomorphism) $\Om_{j,k}$ for $j\in\{1,\dots,n_k\}$ and $k\in\{1,\dots,n\}$, where some of the domains $\Om_{1,k}$ may agree, satisfying:
\begin{itemize}
\item[$i)$] each polyhedron $\Om_{j,k}$ is a $\nabla u_{j,k}$-admissible domain with $u_{j,k}\in C^2(\overline{\Om_{j,k}})$;
\item[$ii)$] each internal face of the chain $\Om_{1,k}\to\Om_{2,k}\to\dots\to\Om_{n_k,k}$ made of $n_k$ contiguous domains, is an inflow boundary for one domain and an outflow boundary for the contiguous domain, or equivalently
\beq\label{inf-outf}
{\partial u_{j,k}\over\partial\nu}\,{\partial u_{j-1,k}\over\partial\nu}>0\quad\mbox{on }\partial\Om_{j,k}\cap\partial\Om_{j-1,k}
\quad\mbox{for any }j\in\{2,\dots,n_k\},
\eeq
where $\nu$ is the outer normal of $\partial\Om_{j,k}$;
\item[$iii)$] each external face of the chain $\Om_{1,k}\to\Om_{2,k}\to\dots\to\Om_{n_k,k}$ is
\begin{itemize}
\item either a boundary part of~$\partial\Om$,
\item or a surface tangential to some $\nabla u_{j,k}$,
\item or an inflow or outflow boundary of $\Om_{1,k}$ which is (possibly) connected to another chain $\Om_{1,k}=\Om_{1,j}\to\Om_{2,j}\to\dots\to\Om_{n_j,j}$.
\end{itemize}
\end{itemize}
\end{Def}
\begin{Exa}
\hfill
\begin{enumerate}
\item Figure~\ref{fig2} represents an admissible domain $\Om$ composed of the $n=4$ chains
\[
\left\{\ba{l}
\Om_{1,1}\to\Om_{2,1}\to\Om_{3,1}\to\Om_{4,1}
\\
\Om_{1,1}=\Om_{1,2}\to\Om_{2,2}
\\
\Om_{1,1}=\Om_{1,3}\to\Om_{2,3}\to\Om_{3,3}
\\
\Om_{1,4}\to\Om_{2,4}.
\ea\right.
\]
The three first chains are connected to the same set $\Om_{1,1}$. The fourth one is separated from three others by surfaces which are tangential to the gradient.
\item The domain $\Om$ of Figure~\ref{fig3} is composed of $n=1$ chain made of $4$ $\nabla u_k$-admissible sets. It is not admissible, since the chain $\Om_1\to\Om_2\to\Om_3\to\Om_4$ has an external boundary which is neither a boundary part of~$\partial\Om$ nor a surface tangential to some gradient $\nabla u_{k}$. This creates a conflict for defining a suitable conductivity $\si_k$ in each domain $\Om_k$ (see Remark~\ref{rem.fig23}, 2. below).
\end{enumerate}
\end{Exa}
\begin{Thm}\label{thm.jump}
Let $\Om$ be an admissible domain composed of $\nabla u_{j,k}$-admissible open sets $\Om_{j,k}$ for $j\in\{1,\dots,n_k\}$ and $k\in\{1,\dots,n\}$, according to Definition~\ref{def.adm}, and let $u\in C(\overline{\Om})$ be such that $u=u_{j,k}$ in $\overline{\Om_{j,k}}$.
Then, there exists a piecewise continuous positive conductivity $\si$ such that
\beq\label{siOmkj}
\left\{\ba{ll}
\si_{\mid\overline{\Om_{k,j}}}\in C^1\big(\overline{\Om_{k,j}}\big) & \mbox{for }j\in\{1,\dots,n_k\}\mbox{ and }k\in\{1,\dots,n\},
\\ \ecart
{\rm div}\left(\si\nabla u\right)=0 & \mbox{in }\D'(\Om).
\ea\right.
\eeq
\par
Conversely, let $\Om$ be a bounded domain of $\RR^d$ composed of $n$ generalized polyhedra $\Om_k$, and let $u$ be a function in $C(\overline{\Om})$ such that $u_k:=u_{\mid\overline{\Om_k}}\in C^2(\overline{\Om_k})$ and $\Om_k$ is a $\nabla u_k$-admissible domain for $k\in\{1,\dots,n\}$. Assume that $\si$ is a positive function in $C(\overline{\Om})$ such that $\si_k:=\si_{\mid\overline{\Om_k}}\in C^1(\overline{\Om_k})$ and ${\rm div}\,\left(\si\nabla u\right)=0$ in $\D'(\Om)$. Then, for any contiguous polyhedra $\Om_j$ and $\Om_k$, the common face $\Ga_{j,k}:=\partial\Om_j\cap\partial\Om_k$ is either a surface tangential to $\nabla u$, or an inflow (resp. outflow) boundary of $\Om_j$ and an outflow (resp. inflow) boundary of $\Om_k$.
\end{Thm}

\begin{figure}
\centering
\vskip -0.cm
\includegraphics[scale=.5]{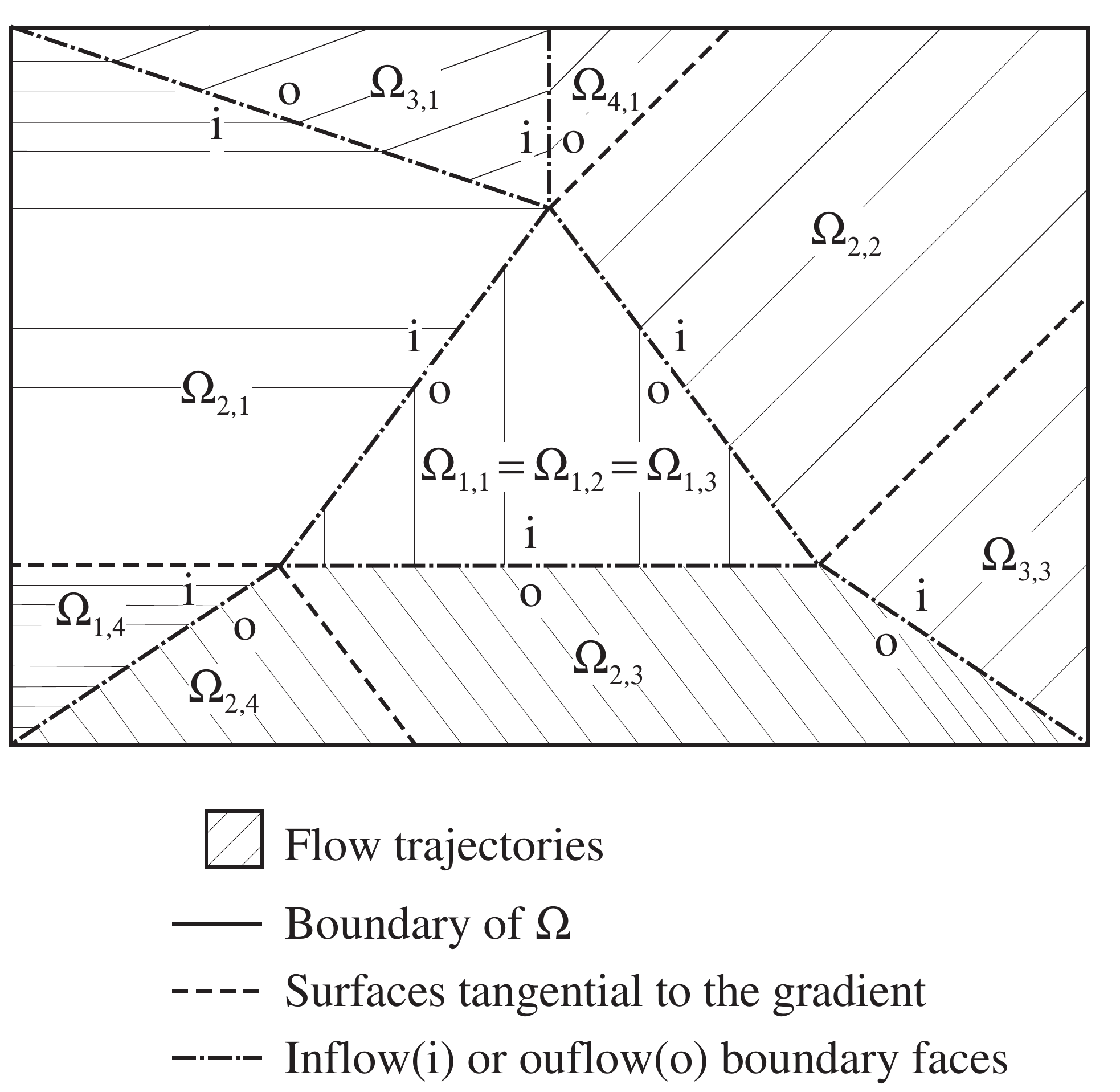}
\vskip 0.cm
\caption{\it An admissible domain $\Om$ composed of $n=4$ chains}
\label{fig2}
\end{figure}

\begin{figure}
\centering
\vskip -0.cm
\includegraphics[scale=.5]{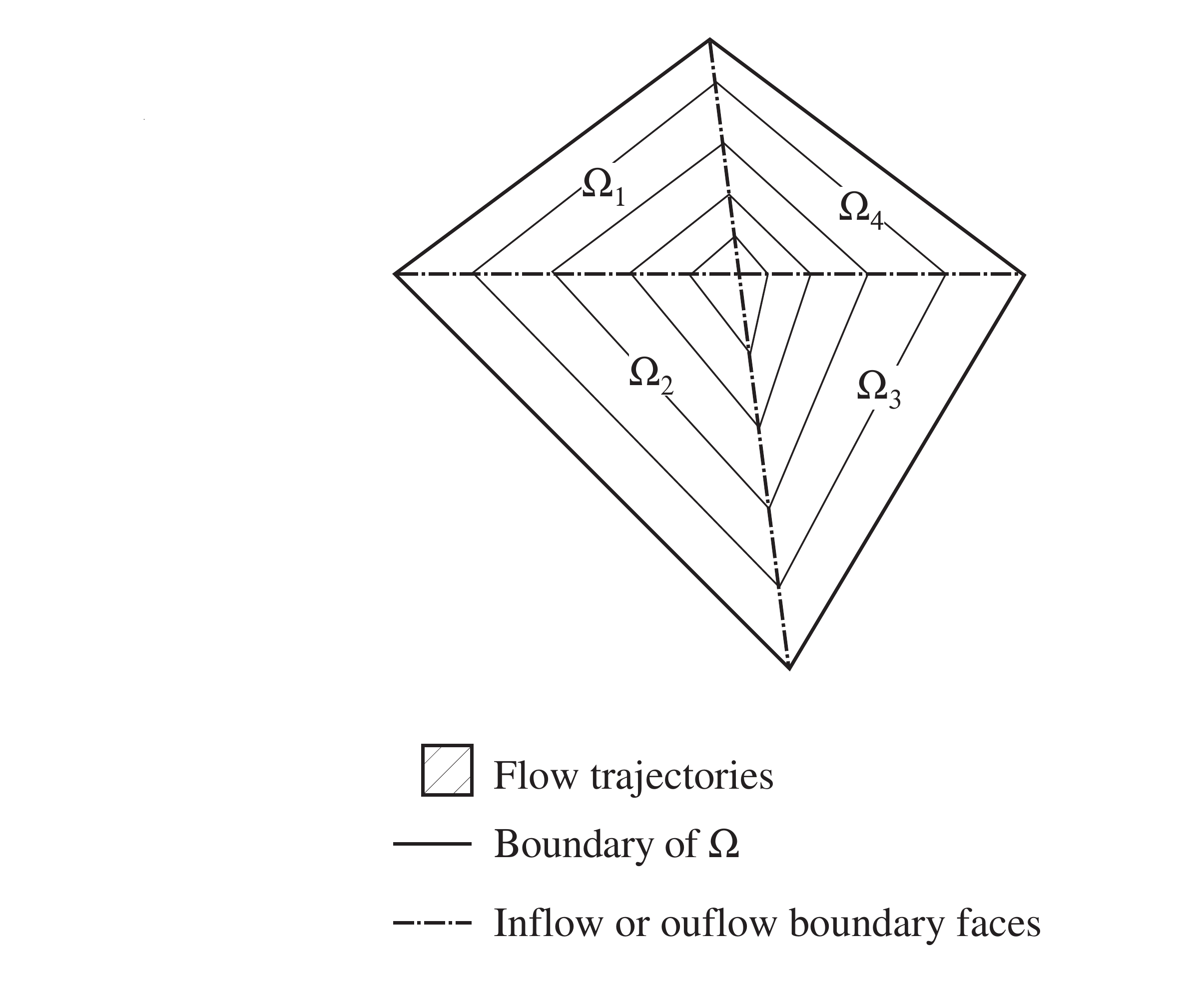}
\vskip 0.cm
\caption{\it A non-admissible domain $\Om$ with $n=1$ chain: $\Om_1\to\Om_2\to\Om_3\to\Om_4$}
\label{fig3}
\end{figure}

\noindent
{\bf Proof of Theorem~\ref{thm.jump}.}
The idea is to construct in each chain $\Om_{1,k}\to\Om_{2,k}\to\dots\to\Om_{n_k,k}$ for $k\in\{1,\dots,n\}$, successively the conductivities $\si_{1,k},\dots,\si_{n_k,k}$. To this end, the conductivity $\si_{j-1,k}$ being constructed in the domain $\Om_{j-1,k}$ for some $j\in\{2,\dots,n_k\}$, we will choose a suitable positive continuous function $\ga_{j,k}$ on the inflow or outflow boundary face $\partial\Om_{j,k}\cap\partial\Om_{j-1,k}$, which
\begin{itemize}
\item determines the conductivity $\si_{j,k}$ in the $\nabla u_{j,k}$-admissible domain $\Om_{j,k}$ by Proposition~\ref{pro.BVR},
\item satisfies the flux continuity condition through the surface $\partial\Om_{j,k}\cap\partial\Om_{j-1,k}$.
\end{itemize}
\par
For $k\in\{1,\dots,n\}$, fix the conductivity equal to $1$ on the inflow or outflow boundary face of~$\Om_{1,k}$, which by Proposition~\ref{pro.BVR} determines a unique conductivity $\si_{1,k}\in C^1(\overline{\Om_{1,k}})$ such that ${\rm div}\left(\si_{1,k}\nabla u\right)=0$ in~$\Om_{1,k}$.
\par
Next, using an induction argument we will construct a suitable piecewise continuous conductivity along the chain $\Om_{1,k}\to\cdots\to\Om_{n_k,k}$. Assume that for some $j\in\{2,\dots,n_k\}$, we have built a piecewise conductivity $\si=\si_{i,k}$ in $\overline{\Om_{i,k}}$ for $i\in\{1,\dots,j-1\}$, solution to the equation
\[
{\rm div}\left(\si\nabla u\right)=0\;\;\mbox{in }{\rm int}\!\left[\Om_{1,k}\cup\,\cup_{i=2}^{j-1}\;(\Om_{i,k}\cup \Ga_{i,k})\right],
\]
where $\Ga_{i,k}:=\partial\Om_{i,k}\cap\partial\Om_{i-1,k}$ is the common face of $\Om_{j,k}$ and $\Om_{j-1,k}$.
By the condition \eqref{inf-outf} on~$\Ga_{j,k}$ there exists a positive function $\ga_{j,k}\in C(\Ga_{j,k})$ such that
\beq\label{Gajk}
\ga_{j,k}\,{\partial u_{j,k}\over\partial\nu}=\si_{j-1,k}{\partial u_{j-1,k}\over\partial\nu}\quad\mbox{on }\Ga_{j,k},
\eeq
where $\nu$ is the outer normal of $\partial\Om_{j,k}$. Since by the assumption $ii)$ of Definition~\ref{def.adm} $\Ga_{j,k}$ is an inflow or outflow boundary face of the $\nabla u_{j,k}$-admissible domain $\Om_{j,k}$, by Proposition~\ref{pro.BVR}  there exists a positive conductivity $\si_{j,k}\in C(\overline{\Om_{j,k}})$ taking the value $\ga_{j,k}$ on $\Ga_{j,k}$ and solution to the equation ${\rm div}\left(\si_{j,k}\nabla u\right)=0$ in $\Om_{j,k}$. Then, equality \eqref{Gajk} reads as the flux continuity condition through $\Ga_{j,k}$.
It follows that the conductivity $\si:=\si_{i,k}$ in $\overline{\Om_{i,k}}$ for $i\in\{1,\dots,j\}$, is solution to the equation
\[
{\rm div}\left(\si\nabla u\right)=0\;\;\mbox{in }{\rm int}\!\left[\Om_{1,k}\cup\,\cup_{i=2}^{j}\;(\Om_{i,k}\cup \Ga_{i,k})\right],
\]
which concludes the induction proof.
Therefore, we has just constructed a piecewise continuous positive function
\beq\label{piecewisesi}
\si=\si_{j,k}\mbox{ in }\overline{\Om_{j,k}}\quad\mbox{solution to}\quad
{\rm div}\left(\si\nabla u\right)=0\;\;\mbox{in }
{\rm int}\!\left(\Om_{1,k}\cup\,\cup_{j=2}^{n_k}\;(\Om_{j,k}\cup \Ga_{j,k})\right).
\eeq
\par
Now, according to Definition~\ref{def.adm} consider the partition $(K_i)_{1\leq i\leq p}$ of $\{1,\dots,n\}$ such that the sets $\Om_{1,k}$ agree to the same set $\Om_{1,k_i}$ ($k_i\in K_i$) for any $k\in K_i$ and $i\in\{1,\dots,p\}$.
Since for each $i\in\{1,\dots,p\}$ the chains $\Om_{1,k}\to\Om_{2,k}\to\cdots\to\Om_{n_k,k}$ are connected to the set $\Om_{1,k_i}$ for any $k\in K_i$, by the definition~\eqref{piecewisesi} of the piecewise continuous conductivity $\si$ we thus have
\beq\label{divsiDu=01}
{\rm div}\left(\si\nabla u\right)=0\quad\mbox{in }{\rm int}\!\left(\bigcup_{k\in K_i}\left[\Om_{1,k_i}\cup\,\cup_{j=2}^{n_k}\;(\Om_{j,k}\cup \Ga_{j,k})\right]\right)
\quad\mbox{for any }i\in\{1,\dots,p\}.
\eeq
Moreover, by the assumption $iii)$ of Definition~\ref{def.adm} we have
\beq\label{siDunu=02}
{\partial u\over\partial\nu}=0\quad\mbox{on }\partial\left(\bigcup_{k\in K_i}\left[\Om_{1,k_i}\cup\,\cup_{j=2}^{n_k}\;(\Om_{j,k}\cup \Ga_{j,k})\right]\right)\setminus\partial\Om
\quad\mbox{for any }i\in\{1,\dots,p\}.
\eeq
Let $\ph\in C^\infty_c(\Om)$. Therefore, integrating by parts and using \eqref{divsiDu=01}, \eqref{siDunu=02} we get that
\[
\int_\Om\si\nabla u\cdot\nabla\ph\,dx=\sum_{i=1}^p\int_{\bigcup_{k\in K_i}\left[\Om_{1,k_i}\cup\,\cup_{j=2}^{n_k}\;(\Om_{j,k}\cup \Ga_{j,k})\right]}\si\nabla u\cdot\nabla\ph\,dx=0,
\]
which implies that the piecewise continuous conductivity $\si$ of \eqref{piecewisesi} is solution to the equation ${\rm div}\left(\si\nabla u\right)=0$ in $\D'(\Om)$.
\par\ms
Conversely, let $\Om$ be a bounded domain of $\RR^d$ composed of $n$ generalized polyhedra $\Om_k$ for $k\in\{1,\dots,n\}$. Let  $u\in C(\overline{\Om})$ be such that $u_k:=u_{\mid\overline{\Om_k}}\in C^2(\overline{\Om_k})$, and $\Om_k$ is $\nabla u_k$-admissible. Assume that $\si$ is a positive piecewise continuous function such that $\si_k:=\si_{\mid\overline{\Om_k}}\in C^1(\overline{\Om_k})$ and ${\rm div}\,\left(\si\nabla u\right)=0$ in $\D'(\Om)$. Consider two contiguous polyhedra $\Om_j$ and $\Om_k$, the common face of which $\Ga_{j,k}:=\partial\Om_j\cap\partial\Om_k$ is not a surface tangential to $\nabla u$. The flux continuity condition through $\Ga_{j,k}$ reads as
\beq\label{fluxcontGajk}
\si_j\,{\partial u_j\over\partial\nu}=\si_k\,{\partial u_k\over\partial\nu}\quad\mbox{on }\Ga_{j,k},
\eeq
where $\nu$ is the outer normal to $\partial\Om_j$, which implies that
\[
{\partial u_j\over\partial\nu}\,{\partial u_k\over\partial\nu}>0\quad\mbox{on }\Ga_{j,k}.
\]
Therefore, $\Ga_{j,k}$ is an inflow (resp. outflow) boundary face of $\Om_j$, and an outflow (resp. inflow) boundary face of $\Om_k$.
The proof of Theorem~\ref{thm.jump} is now complete.
\cqfd
\begin{Rem}\label{rem.fig23}
\hfill
\begin{enumerate}
\item In the case of Figure~\ref{fig2} the domain $\Om$ is composed of $9$ polyhedra $\Om_{j,k}$ grouped into $4$ chains with $11$ internal faces.
The step by step construction of Theorem~\ref{thm.jump} reads as follows:
\begin{itemize}
\item We prescribe the conductivity on the say inflow face $\partial\Om_{1,1}\cap\partial\Om_{2,3}$ of $\Om_{1,1}$, which determines the conductivity $\si_{1,1}$. Then, $\partial\Om_{1,1}\cap\partial\Om_{2,1}$ and $\partial\Om_{1,1}\cap\partial\Om_{2,2}$ are outflow faces of~$\Om_{1,1}$.
\item We choose successively the conductivities on the inflow face $\partial\Om_{1,1}\cap\partial\Om_{2,1}$ of $\Om_{2,1}$, the outflow face $\partial\Om_{2,1}\cap\partial\Om_{3,1}$ of $\Om_{3,1}$, and the outflow face $\partial\Om_{3,1}\cap\partial\Om_{4,1}$ of $\Om_{4,1}$, which determine the conductivities $\si_{2,1},\si_{3,1},\si_{4,1}$ ensuring the flux continuity conditions on $\partial\Om_{1,1}\cap\partial\Om_{2,1}$, $\partial\Om_{2,1}\cap\partial\Om_{3,1}$, $\partial\Om_{3,1}\cap\partial\Om_{4,1}$.
\item We choose the conductivity on the inflow face $\partial\Om_{1,1}\cap\partial\Om_{2,2}$ of $\Om_{2,2}$, which determines the conductivity $\si_{2,2}$ ensuring the flux continuity condition on $\partial\Om_{1,1}\cap\partial\Om_{2,2}$.
\item We choose successively the conductivities on the outflow face $\partial\Om_{1,1}\cap\partial\Om_{2,3}$ of $\Om_{2,3}$ and the inflow face $\partial\Om_{2,3}\cap\partial\Om_{3,3}$ of $\Om_{3,3}$, which determine the conductivities $\si_{2,3},\si_{3,3}$ ensuring the flux continuity conditions on $\partial\Om_{1,1}\cap\partial\Om_{2,3}$, $\partial\Om_{2,3}\cap\partial\Om_{3,3}$.
\item We prescribe the conductivity on the say inflow face $\partial\Om_{1,4}\cap\partial\Om_{2,4}$ of $\Om_{1,4}$, which determines the conductivity $\si_{1,4}$. Then, we choose the conductivity on the ouflow face $\partial\Om_{1,4}\cap\partial\Om_{2,4}$ of $\Om_{2,4}$, which determines the conductivity $\si_{2,4}$ ensuring the flux continuity condition on $\partial\Om_{1,4}\cap\partial\Om_{2,4}$.
\item The $4$ remaining faces $\partial\Om_{4,1}\cap\partial\Om_{2,2}$, $\partial\Om_{2,2}\cap\partial\Om_{3,3}$, $\partial\Om_{2,3}\cap\partial\Om_{2,4}$, $\partial\Om_{2,1}\cap\partial\Om_{1,4}$ are tangential to the gradient, and thus satisfy the flux continuity conditions.
\end{itemize}
\item In the case of Figure~\ref{fig2} the domain $\Om$ is made of one chain composed of $4$ polyhedra. For example, we prescribe the conductivity on the say inflow face $\partial\Om_1\cap\partial\Om_2$ of $\Om_1$. Then, the flux continuity conditions on the faces $\partial\Om_1\cap\partial\Om_2$, $\partial\Om_2\cap\partial\Om_3$, $\partial\Om_3\cap\partial\Om_4$ determine successively the conductivities $\si_k$ in $\Om_k$ for $k=1,2,3,4$. But then the flux continuity condition on the face $\partial\Om_1\cap\partial\Om_4$ does not hold in general.
\end{enumerate}
\end{Rem}
}
\section{Examples}\label{s.exa}
\subsection{Example 1}
Let $\Om$ be an open set of $\RR^2$ which is star-shaped with respect to the origin. Let $\xi_1,\dots,\xi_n$ be $n\geq 2$ non-zero vectors of $\RR^2$ such that the open cones
\beq\label{Omkxik}
\left\{\ba{ll}
\Om_k:=\big\{s\,\xi_k+t\,\xi_{k+1},\,,s,t>0\big\} & \mbox{for }1\leq k\leq n-1
\\ \ecart
\Om_n:=\big\{s\,\xi_1+t\,\xi_n,\,,s,t>0\big\} & \mbox{for }k=n,
\ea\right.
\eeq
do not contain any vector $\xi_j$.
\begin{figure}[t!]
\centering
\vskip -0.cm
\includegraphics[scale=.5]{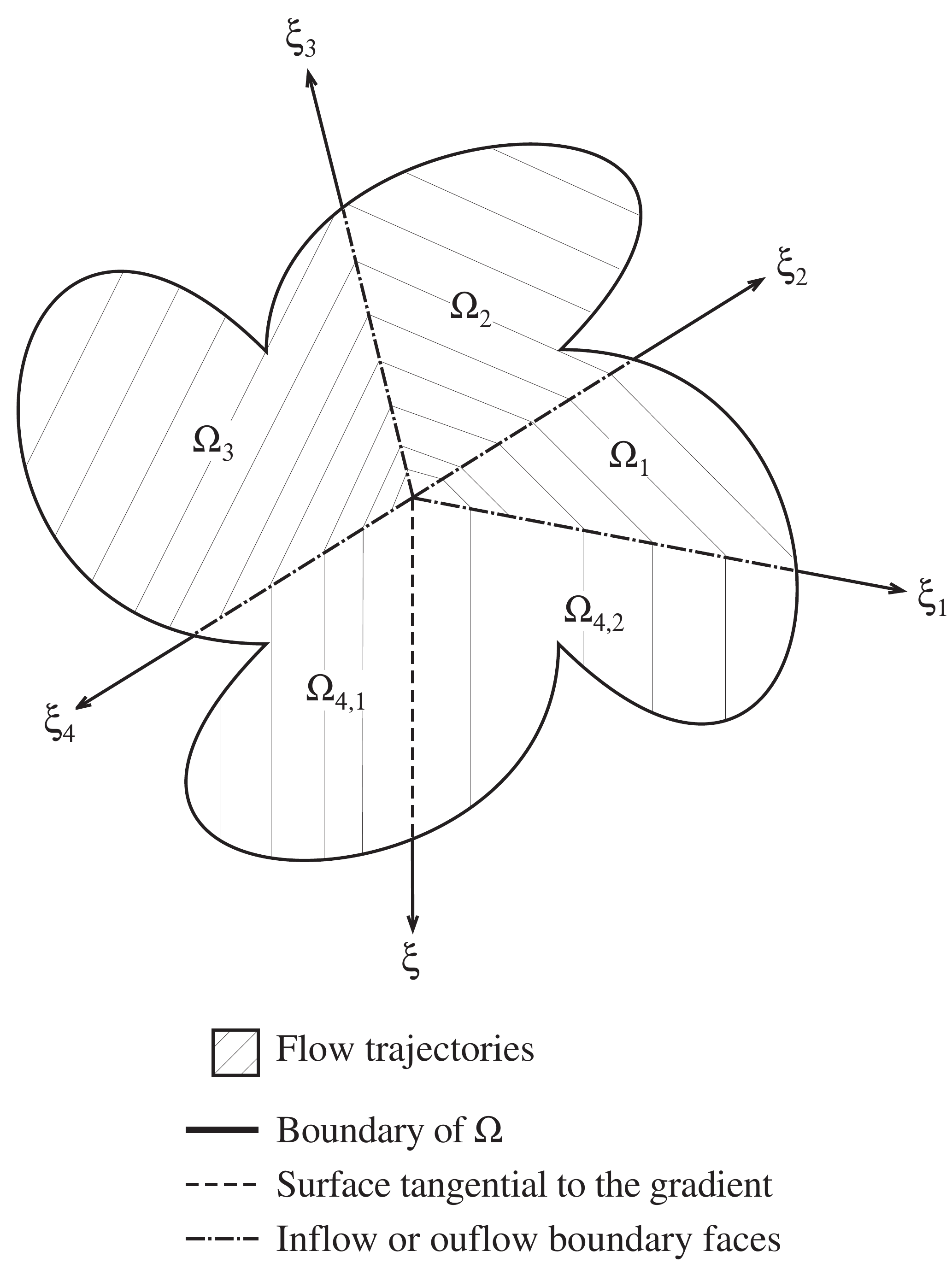}
\vskip 0.cm
\caption{\it Triangulation of $\Om$ by the cones $\Om_1$, $\Om_2$, $\Om_3$, and $\Om_4={\rm int}\left(\overline{\Om_{4,1}}\cup\overline{\Om_{4,2}}\right)$ with $\xi\parallel\la_4$}
\label{fig4}
\end{figure}

Consider a function $u\in C(\overline{\Om})$ of finite element type $\PP_1$ (see, {\em e.g.} \cite[Section~2.2]{Cia}, {\em i.e.} there exists constant vectors $\la_k\in\RR^2$ such that
\beq\label{Dulak}
\nabla u=\la_k\;\;\mbox{in }\Om_k\quad\mbox{for }k\in\{1,\dots,n\}.
\eeq
This imposes the flux continuity conditions
\beq\label{conlak1}
(\la_k-\la_{k-1})\cdot\xi_k=0,\;\;\forall\,k\in\{2,\dots,n\}\quad\mbox{and}\quad(\la_1-\la_n)\cdot\xi_1=0.
\eeq
Up to decrease the value of $n$ we can also assume that
\beq\label{conlak2}
\la_k-\la_{k-1}\neq 0,\;\;\forall\,k\in\{2,\dots,n\}\quad\mbox{and}\quad\la_1-\la_n\neq 0.
\eeq
\par
Similarly to the case of Figure~\ref{fig3} (see Remark~\ref{rem.fig23}, 2.) the chain $\Om_{1}\to\Om_{2}\to\dots\to\Om_{n}$ does not satisfy the condition $iii)$ of Definition~\ref{def.adm}.  Indeed, the existence of constant conductivities $\si_k$ in $\Om_k$ satisfying the flux continuity condition~\eqref{fluxcontGajk} reads as
\[
\si_k\det\left(\xi_k,\la_k\right)=\si_{k-1}\det\left(\xi_k,\la_{k-1}\right),\;\;\forall\,k\in\{2,\dots,n\}\quad\mbox{and}\quad
\si_n\det\left(\xi_1,\la_n\right)=\si_1\det\left(\xi_1,\la_1\right),
\]
which thus implies the constraint
\beq\label{sikxiklak}
\prod_{k=1}^n\,\det\left(\xi_k,\la_k\right)=\det\left(\xi_1,\la_n\right)\,\prod_{k=2}^{n}\,\det\left(\xi_k,\la_{k-1}\right).
\eeq
\par
A less restrictive alternative is to assume that for some $k\in\{1,\dots,n\}$, say $k=n$ without loss of generality, there exists a vector $\xi\in\RR^2$ satisfying
\beq\label{xi}
\xi\in\Om_n\setminus\{0\}\quad\mbox{and}\quad\xi\parallel\la_n.
\eeq
Hence, defining the subsets of $\Om_n$
\[
\Om_{n,1}:=\big\{s\,\xi+t\,\xi_n,\,,s,t>0\big\}\quad\mbox{and}\quad\Om_{n,2}:=\big\{s\,\xi+t\,\xi_1,\,,s,t>0\big\},
\]
we have
\beq\label{Dunu=012}
{\partial u\over\partial\nu}=0\quad\mbox{on }\partial\Om_{n,1}\cap\partial\Om_{n,2}\subset\RR\,\xi.
\eeq
Therefore, by \eqref{Dulak} and \eqref{Dunu=012} the chain $\Om_{n,2}\to\Om_1\to\dots\to\Om_{n-1}\to\Om_{n,1}$ satisfies the conditions $i)$ and $iii)$ of Definition~\ref{def.adm} (see Figure~\ref{fig4} and compare to Figure~\ref{fig3}).
Then, taking into account conditions \eqref{conlak1} and \eqref{conlak2} the condition $ii)$ of Definition~\ref{def.adm} is equivalent to
\beq\label{conlak}
\det\left(\xi_k,\la_k\right)\det\left(\xi_k,\la_{k-1}\right)>0,\;\;\forall\,k\in\{2,\dots,n\}\quad\mbox{and}\quad
\det\left(\xi_1,\la_1\right)\det\left(\xi_1,\la_n\right)>0.
\eeq
Therefore, by Theorem~\ref{thm.jump} $\nabla u$ is isotropically realizable in $\Om$ if and only if condition \eqref{conlak} holds true.
Finally, due to condition \eqref{conlak} a suitable piecewise constant conductivity is given by
\beq\label{silai}
\si=\left\{\ba{cll}
\dis {\det\left(\xi_1,\la_n\right)\over\det\left(\xi_1,\la_1\right)} & \mbox{in }\Om_1 &
\\ \ecart
\dis {\det\left(\xi_1,\la_n\right)\over\det\left(\xi_1,\la_1\right)}\,\prod_{j=2}^k{\det\left(\xi_j,\la_{j-1}\right)\over\det\left(\xi_j,\la_j\right)} & \mbox{in }\Om_k & \mbox{for }2\leq k\leq n-1
\\ \ecart
\dis {\det\left(\xi_1,\la_n\right)\over\det\left(\xi_1,\la_1\right)}\,\prod_{j=2}^n{\det\left(\xi_j,\la_{j-1}\right)\over\det\left(\xi_j,\la_j\right)} & \mbox{in }\Om_{n,1} &
\\ \ecart
\dis 1 & \mbox{in }\Om_{n,2}. &
\ea\right.
\eeq

\begin{Rem}\label{rem.exa}
We can also extend the previous two-dimensional example to dimension three replacing the open cones \eqref{Omkxik} as follows.
Let $\Om$ be an open set of $\RR^3$ which is star-shaped with respect to the origin.
Let $\xi_1,\dots,\xi_n$ be $n\geq 3$ non-zero vectors of $\RR^3$ such that the open cones
\beq\label{Omjkl}
\Om_{i,j,k}:=\Om\cap\big\{r\,\xi_i+s\,\xi_j+t\,\xi_k,\,r,s,t>0\big\}\quad\mbox{if }\det\left(\xi_i,\xi_j,\xi_k\right)\neq 0,
\eeq
do not contain any vector $\xi_\ell$. For example, if $(e_1,e_2,e_3)$ is a basis of $\RR^3$ and $n=6$ with
\[
\xi_1=e_1,\ \xi_2=e_2,\ \xi_3=e_3,\ \xi_4=-e_1,\ \xi_5=-e_2,\ \xi_6=-e_3,
\]
there are $8$ open cones of type \eqref{Omjkl}.
\end{Rem}
\subsection{Example 2}
Let $f$ be a function in $W^{2,\infty}_{\rm loc}(\RR^{d-1})$ for $d\geq 2$, and let $g,h$ be $2$ functions in $C^2(\RR)$ such that
\beq\label{fgh}
\left\{\ba{l}
f\mbox{ satisfies condition \eqref{nojumpDu} in }\RR^{d-1},
\\ \ecart
g(0)=h(0),
\\ \ecart
g',h'\mbox{ are uniformly continuous in }\RR\;\;\mbox{and}\;\;g'(t)\,h'(t)\neq 0,\;\;\forall\,t\in\RR.
\ea\right.
\eeq
Consider the function $u\in C(\RR^d)$ defined by
\beq\label{ufgh}
u(x)=\left\{\ba{ll}
u_1(x_1,x'):=g(x_1)+f(x') & \mbox{if }(x_1,x')\in\Om_1:=(0,\infty)\times\RR
\\ \ecart
u_2(x_1,x'):=h(x_1)+f(x') & \mbox{if }(x_1,x')\in\Om_2:=(-\infty,0)\times\RR,
\ea\right.
\eeq
so that $u$ satisfies the conditions $i)$ and $iii)$ (which is empty there) of Definition~\ref{def.adm}. Moreover, the function $\nabla u$ is piecewise continuous in $\RR^d$, and condition $ii)$ of Definition~\ref{def.adm} is reduced to
\beq\label{g'h'0}
g'(0)\,h'(0)>0.
\eeq
Due to the separation of the variables $x_1$ and $x'$, the gradient flow $X=(X_1,X')$ associated with $\nabla u_1$ satisfies
\[
\left\{\ba{rl}
\dis {\partial X_1\over \partial t}(t,x_1) & =g'\big(X_1(t,x_1)\big)
\\ \ecart
X_1(0,x_1) & =x_1,
\\ \ecart
\dis {\partial X'\over \partial t}(t,x') & =\nabla_{x'}f\big(X'(t,x)\big)
\\ \ecart
X'(0,x') & =x'
\ea\right.
\quad\mbox{for }t\in\RR,\ x=(x_1,x')\in\RR^d,
\]
which yields
\beq\label{XDu1}
\left\{\ba{rl}
X_1(t,x_1) & =G^{-1}\big(t+G(x_1)\big)
\\ \ecart
X_1(0,x_1) & =x_1,
\\ \ecart
\dis {\partial X'\over \partial t}(t,x') & =\nabla_{x'}f\big(X'(t,x)\big)
\\ \ecart
X'(0,x) & =x'
\ea\right.
\quad\mbox{for }t\in\RR,\ x=(x_1,x')\in\RR^d,
\eeq
where $G^{-1}$ is the inverse function of the primitive $G$ of $1/g'$ in $\RR$ such that $G(0)=0$.
For a.e. $x\in\RR^d$, the flow $X(\cdot,x)$ reaches the surface $\{x_1=0\}$ at the time $\tau_1(x)=-G(x_1)$ which implies $X_1\big(\tau_1(x),x_1\big)=0$. Then, by Theorem~\ref{thm.nojump} and formula~\eqref{sitau} with $u_1$, for any constant $\la>0$, the gradient $\nabla u_1$ is realizable  with the continuous conductivity
\[
\si_1(x)=\la\exp\left(\int_0^{-G(x_1)}\left[g''\big(X_1(s,x_1)\big)+\De_{x'}f\big(X'(s,x')\big)\right]ds\right)\quad\mbox{for }x\in\RR^d,
\]
which using the change of variable $t=X_1(s,x_1)=G^{-1}\big(s+G(x_1)\big)$ yields
\beq\label{si1}
\si_1(x)=\la\,{g'(0)\over g'(x_1)}\,\exp\left(\int_0^{-G(x_1)}\De_{x'}f\big(X'(s,x')\big)\,ds\right)\quad\mbox{for a.e. }x\in\RR^d.
\eeq
Similarly, the gradient $\nabla u_2$ is realizable in $\RR^d$ with the continuous conductivity
\beq\label{si2}
\si_2(x)={h'(0)\over h'(x_1)}\,\exp\left(\int_0^{-H(x_1)}\De_{x'}f\big(X'(s,x')\big)\,ds\right)\quad\mbox{for a.e. }x\in\RR^d,
\eeq
where $H$ is the primitive of $1/h'$ in $\RR$ such that $H(0)=0$.
Choosing $\la=h'(0)/g'(0)>0$ by \eqref{g'h'0}, we get the flux continuity condition across the interface $\{x_1=0\}$, {\em i.e.}
\[
\si_1(0,x')\,{\partial u_1\over\partial x_1}(0,x')=\si_2(0,x')\,{\partial u_2\over\partial x_1}(0,x')=h'(0)\quad\mbox{for }x'\in\RR^{d-1}.
\]
Therefore, the gradient $\nabla u$ is realizable with the piecewise continuous conductivity
\beq\label{sisi1si2}
\si(x)=\left\{\ba{ll}
\dis {h'(0)\over g'(x_1)}\,\exp\left(\int_0^{-G(x_1)}\De_{x'}f\big(X'(s,x')\big)\,ds\right) & \mbox{if }x\in(0,\infty)\times\RR^{d-1}
\\ \ecart
\dis {h'(0)\over h'(x_1)}\,\exp\left(\int_0^{-H(x_1)}\De_{x'}f\big(X'(s,x')\big)\,ds\right) & \mbox{if }x\in(-\infty,0)\times\RR^{d-1}.
\ea\right.
\eeq
\par\bs\noindent
{\bf Acknowledgements.}
The author is very grateful to the unknown referees for their careful reading and quite relevant comments and references, which have significantly improved the presentation of the paper especially for Section~\ref{s.jump}.
\end{document}